\documentclass{article}%
\usepackage{amsmath}
\usepackage{makeidx}%
\usepackage{amsfonts}%
\usepackage{amssymb}%
\usepackage{graphicx}
\providecommand{\U}[1]{\protect\rule{.1in}{.1in}}
\newtheorem{theorem}{Theorem}

\newtheorem{corollary}[theorem]{Corollary}

\newtheorem{lemma}[theorem]{Lemma}

\begin{document}

\title{\textbf{Inverses of Motzkin and Schr\"{o}der Paths}}
\author{Heinrich Niederhausen\\Florida Atlantic University, Boca Raton}
\date{May 2011}
\maketitle

\begin{abstract}
We suggest three applications for the inverses: For the inverse Motzkin matrix
we look at Hankel determinants, and counting the paths inside a horizontal
band, and for the inverse Schr\"{o}der matrix we look at the paths inside the
same band, but ending on the top side of the band.

\end{abstract}

\section{Introduction}

We adopt the convention that lattice paths without restrictions are called
\textquotedblleft\emph{Grand}\textquotedblright; the \emph{Grand Catalan
numbers} (step set $\left\{  \nearrow,\searrow\right\}  $) are the number of
paths from the origin, taking only $\nearrow$ and $\searrow$ steps, and ending
on the $x$-axis at $\left(  2n,0\right)  $. The Grand Catalan numbers are the
\emph{Central Binomial coefficients}, $\binom{2n}{n}$, with generating
function $1/\sqrt{1-4t^{2}}=\sum_{n\geq0}\binom{2n}{n}t^{2n}$. The
\emph{wheighted} \emph{Grand Motzkin} numbers $G_{n}$ take steps from
$\left\{  \nearrow,\searrow,\longrightarrow\right\}  $, and end on the
$x$-axis in $\left(  n,0\right)  $. The horizontal steps get the weight
$\omega$. Their generating function is
\begin{equation}
g\left(  t\right)  :=\sum_{n\geq0}G_{n}t^{n}=1/\sqrt{\left(  1-\omega
t\right)  ^{2}-4t^{2}},\label{(GFgrandMotzkin)}%
\end{equation}
and it is seen immediately that for $\omega=0$ the Grand Catalan numbers are
recovered. If $\omega=2$, the $1/\sqrt{\left(  1-2t\right)  ^{2}-4t^{2}%
}=\allowbreak1/\sqrt{1-4t}$ is again a generating function for the Grand
Catalan numbers, but we get $\sum_{n\geq0}\binom{2n}{n}t^{n}$. The
\emph{general} Grand Motzkin numbers $G\left(  n,j\right)  $ enumerate all
paths to $\left(  n,j\right)  $, and the first few are given in the following table.%

\[%
\begin{tabular}
[c]{l||llllll}%
$\uparrow j$ &  &  &  &  &  & \\
7 &  &  &  &  &  & \\
6 &  &  &  &  &  & \\
5 &  &  &  &  &  & $1$\\
4 &  &  &  &  & $1$ & $5\omega$\\
3 &  &  &  & $1$ & $4\omega$ & $5+10\omega^{2}$\\
2 &  &  & $1$ & $3\omega_{\searrow}$ & $4+6\omega^{2}$ & $20\omega
+10\omega^{3}$\\
1 &  & $1$ & $2\omega$ & $3+3\omega^{2}\rightarrow$ & $12\omega+4\omega^{3}$ &
$10+30\omega^{2}+5\omega^{4}$\\
0 & $1$ & $\omega$ & $2+\omega^{2}$ & $6\omega+\omega^{3}\nearrow$ &
$6+12\omega^{2}+\omega^{4}$ & $30\omega+20\omega^{3}+\omega^{5}$\\
-1 &  & $1$ & $2\omega$ & $3+3\omega^{2}$ & $12\omega+4\omega^{3}$ &
$10+30\omega^{2}+5\omega^{4}$\\
-2 &  &  & $1$ & $3\omega$ & $4+6\omega^{2}$ & $20+10\omega^{3}$\\
-3 &  &  &  & $1$ & $4\omega$ & $5+10\omega^{2}$\\
-4 &  &  &  &  & $1$ & $5\omega$\\
-5 &  &  &  &  &  & $1$\\
-6 &  &  &  &  &  & \\\hline\hline
$n\rightarrow$ & 0 & 1 & 2 & 3 & 4 & 5\\
\multicolumn{7}{c}{The general Grand Motzkin numbers ($G_{n}$ is given in row
$0$)}%
\end{tabular}
\]

The lower half of the table is the mirror image of the top half; if we write
the table in matrix form, $G\left(  n,j\right)  $ stands in row $n$ and column
$j$, and we obtain a \emph{Riordan matrix} $G$, because $G\left(
n+1,j+1\right)  =G\left(  n,j\right)  +\omega G\left(  n,j+1\right)  +G\left(
n,j+2\right)  $ (see Rogers \cite{Rogers78}, and \cite{Merlini-Rogers97}). It
follows that
\begin{align*}
\sum_{n\geq j}G\left(  n,j\right)  t^{n}  & =\frac{1}{\sqrt{\left(  1-\omega
t\right)  ^{2}-4t^{2}}}\left(  \frac{1}{2t}\left(  1-t\omega-\sqrt{\left(
\omega t-1\right)  ^{2}-4t^{2}}\right)  \right)  ^{j}\\
& =g\left(  t\right)  \left(  \frac{1}{2t}\left(  1-\omega t-1/g\left(
t\right)  \right)  \right)  ^{j}%
\end{align*}

\[%
\begin{tabular}
[c]{l||llll}%
$\uparrow n$ &  &  &  & \\
0 & $1$ &  &  & \\
1 & $\omega$ & $1$ &  & \\
2 & $2+\omega^{2}$ & $2\omega$ & $1$ & \\
3 & $6\omega+\omega^{3}$ & $3+3\omega^{2}$ & $3\omega$ & $1$\\
4 & $6+12\omega^{2}+\omega^{4}$ & $12\omega+4\omega^{3}$ & $4+6\omega^{2}$ &
$4\omega$\\
5 & $30\omega+20\omega^{3}+\omega^{5}$ & $10+30\omega^{2}+5\omega^{4}$ &
$20\omega+10\omega^{3}$ & $5+10\omega^{2}$\\
6 & $20+90\omega^{2}+30\omega^{4}+\omega^{6}$ & $60\omega+60\omega^{3}%
+6\omega^{5}$ & $15+60\omega^{2}+15\omega^{4}$ & $30\omega+20\omega^{3}%
$\\\hline\hline
$j\rightarrow$ & 0 & 1 & 2 & 3\\
\multicolumn{1}{c||}{} & \multicolumn{4}{||l}{The Riordan matrix $G=\left(
G\left(  n,j\right)  \right)  _{\substack{n=0,\dots\\j=0,\dots,n}}$ ($G_{n} $
is given in column $0$)}%
\end{tabular}
\]

If we \emph{restrict} the $\left\{  \nearrow,\searrow,\overset{\omega
}{\longrightarrow}\right\}  $-paths to the first quadrant, they become
\emph{Motzkin paths} $M\left(  n,j\right)  $. We will look at the inverse
$\left(  m_{i,j}\right)  $ of the matrix $M$, and find it useful in some
applications (see also A. Ralston and P. Rabinowitz, 1978 \cite[p.
256]{Ralston}). Especially, the \emph{bounded} Motzkin numbers $M_{n;w}%
^{\left(  k\right)  }$, the number of Motzkin paths staying strictly below the
parallel to the $x$-axis at height $k$, have a generating function expressed
by the inverse $\left(  m_{i,j}\right)  $, through the \emph{inverse Motzkin
polynomial} $m_{k}\left(  t\right)  =\sum_{i=0}^{k}m_{k,i}t^{k-i}$,%
\[
\sum_{n\geq0}M_{n;\omega}^{\left(  k\right)  }t^{n}=\frac{m_{k-1}(t)}%
{m_{k}\left(  t\right)  }%
\]
(see (\ref{(BoundedM)}). That makes us wonder if paths with different lengths
of the horizontal steps $\left(  w,0\right)  $ have similar properties. In the
case of $w=2$ (\emph{Schr\"{o}der paths}) and $\omega=1$ we have a result,
$\mathcal{S}^{\left(  k\right)  }\left(  t\right)  :=$%
\[
\sum_{n\geq0}S_{n}^{\left(  k\right)  }t^{n}=\frac{\left(  1-t\right)
\sum_{i=0}^{\left(  k-2\right)  /2}t^{2i}\left(  -1\right)  ^{i}%
s_{k-2-2i}\left(  t\right)  +\left(  k\operatorname{mod}2\right)  \left(
-1\right)  ^{\left(  k-1\right)  /2}t^{k-1}}{\left(  1-t\right)  \sum
_{i=0}^{\left(  k-1\right)  /2}t^{2i}\left(  -1\right)  ^{i}s_{k-1-2i}\left(
t\right)  +\left(  \left(  k-1\right)  \operatorname{mod}2\right)  \left(
-1\right)  ^{k/2}t^{k}}%
\]
where the Motzkin terms ($M$ and $m$) are replaced by the corresponding
Schr\"{o}der terms ($S$ and $s$), and $s_{i}\left(  t\right)  $ is the
\emph{inverse Schr\"{o}der polynomial}. Perhaps more interesting is the
generating function identity described in Theorem \ref{ThSchroeder},%
\[
t^{-k}\mathcal{S}^{\left(  k\right)  }\left(  t\right)  s_{k-1}\left(
t\right)  =t^{-k}\mathcal{S}^{\left(  k\right)  }\left(  t,k-1\right)
\]
(as power series) where $S^{\left(  k\right)  }\left(  t,k-1\right)  $ is the
generating function of the bounded Schr\"{o}der number ending on $y=k-1$, just
below the upper boundary.

\section{Motzkin Numbers}

Leaving the Grand Motzkin numbers behind, we introduce the restriction of
counting only paths that do not go below the $x$-axis. A \emph{general
weighted Motzkin path} is counted by the recursion
\[
M\left(  n,m;\omega\right)  =M\left(  n-1,m+1;\omega\right)  +\omega M\left(
n-1,m;\omega\right)  +M\left(  n-1,m-1;\omega\right)
\]
for $m\geq0$, and $M\left(  n,m;\omega\right)  =0$ if $m<0$. The numbers
$M\left(  n,m;\omega\right)  $ are weighted counts of all such path from
$\left(  0,0\right)  $ to $\left(  n,m\right)  $, and we give the special name
$M_{n;\omega}$ to the Motzkin numbers $M\left(  n,0;\omega\right)  $. These
numbers (with weight $\omega=1$) have been studied by Th. Motzkin in 1946
\cite{Motzkin}.%

\[%
\begin{tabular}
[c]{l||llllll}%
$\uparrow m$ &  &  &  &  &  & \\
7 &  &  &  &  &  & \\
6 &  &  &  &  &  & \\
5 &  &  &  &  &  & $1$\\
4 &  &  &  &  & $1$ & $5\omega$\\
3 &  &  &  & $1$ & $4\omega$ & $4+10\omega^{2}$\\
2 &  &  & $1$ & $3\omega\;\;\;\;\;_{\searrow}$ & $3+6\omega^{2}$ &
$15\omega+10\omega^{3}$\\
1 &  & $1$ & $2\omega$ & $2+3\omega^{2}\rightarrow$ & $8\omega+4\omega^{3}$ &
$5+20\omega^{2}+5\omega^{4}$\\
0 & $1$ & $\omega$ & $1+\omega^{2}$ & $3\omega+\omega^{3\nearrow}$ &
$2+6\omega^{2}+\omega^{4}$ & $10\omega+10\omega^{3}+\omega^{5}$\\\hline\hline
& 0 & 1 & 2 & 3 & 4 & 5\\
\multicolumn{1}{c||}{$M_{n;\omega}$ is given in row $0.$} &
\multicolumn{1}{||l}{} & \multicolumn{1}{l}{} & \multicolumn{1}{l}{} &
\multicolumn{1}{l}{} & \multicolumn{1}{l}{} & \multicolumn{1}{l}{}%
\end{tabular}
\]

The above table shows that for $\omega=1$ the original Motzkin numbers are
$1,1,2,4,9,21,51,127,\dots$ (sequence A001006 in the On-Line Encyclopedia of
Integer Sequences (OEIS)).

It is well-known that the general $\omega$-weighted Motzkin numbers have the
generating function
\[
\mu\left(  t;j,\omega\right)  :=\sum_{n\geq0}M\left(  n+j,j;\omega\right)
t^{n}=\left(  \frac{1-\omega t-\sqrt{\left(  1-\omega t\right)  ^{2}-4t^{2}}%
}{2t^{2}}\right)  ^{j+1}%
\]
thus%
\begin{equation}
\mu\left(  t\right)  :=\sum_{n\geq0}M_{n;\omega}t^{n}=\sum_{n\geq0}M\left(
n,0;\omega\right)  t^{n}=\frac{1-\omega t-\sqrt{\left(  1-\omega t\right)
^{2}-4t^{2}}}{2t^{2}}\label{(MotzkinGenFun)}%
\end{equation}
is the generating function of the Motzkin numbers, satisfying the quadratic
equation \cite{Bernhart}
\begin{equation}
\mu\left(  t\right)  =1+\omega t\mu\left(  t\right)  +t^{2}\mu\left(
t\right)  ^{2}\label{(quadratic)}%
\end{equation}

Hence%
\[
M_{n+2;\omega}-\omega M_{n+1;\omega}=\sum_{i=0}^{n}M_{i;\omega}M_{n-i;\omega}%
\]
a well-known identity, combinatorially shown by using the \emph{"First Return
Decomposition"}. The generating function (in $t^{2}$) of the Catalan numbers
$C_{n}$ is easily obtained by setting $\omega=0$ in (\ref{(MotzkinGenFun)}),
but it also follows from $\omega=2$%
\[
\frac{1-2t-\sqrt{\left(  1-2t\right)  ^{2}-4t^{2}}}{2t^{2}}=\frac
{1-2t-\sqrt{1-4t}}{2t^{2}}=\sum_{n\geq1}C_{n}t^{n-1}%
\]
(in $t$). Or we can choose $\omega=1$ and get%
\begin{align*}
\left(  1+t\right)  \sum_{n\geq1}C_{n}\left(  \frac{t}{1+t}\right)  ^{n-1}  &
=\sum_{n\geq0}M_{n;1}t^{n}\\
M_{n;1}  & =\sum_{k=0}^{n}\binom{n}{k}\left(  -1\right)  ^{n-k}C_{k+1}%
\end{align*}
For general $\omega$ follows from (\ref{(MotzkinGenFun)}) the explicit
expression%
\[
M_{n;\omega}=\sum_{k=0}^{n/2}\binom{n}{2k}\frac{\omega^{n-2k}}{2k+1}%
\tbinom{2k+1}{k}.
\]

\section{The Inverse}

Define $\phi\left(  t\right)  $ such that $t/\phi\left(  t\right)  $ is the
compositional inverse of $t\mu\left(  t\right)  $ thus%
\[
\phi\left(  t\mu\left(  t\right)  \right)  =\mu\left(  t\right)  =1+\omega
t\mu\left(  t\right)  +t^{2}\mu\left(  t\right)  ^{2}%
\]
by (\ref{(quadratic)}), and therefore%
\[
\phi\left(  t\right)  =1+\omega t+t^{2}%
\]
This simple form of the inverse is the reason for many special results for
Motzkin numbers. Note that
\[
1/\phi\left(  t\right)  =\left(  1+\omega t+t^{2}\right)  ^{-1}=\sum_{n\geq
0}U_{n}\left(  -\omega/2\right)  t^{n}%
\]
the generating function of the \emph{Chebychef polynomials} of the second kind.

Because of the inverse relationship between $t\mu(t)$ and $t/\phi(t)$ we have
that the matrix inverse of $\left(  M\left(  i,j;\omega\right)  \right)
_{n\times n}$ equals $\left(  m_{i,j}\right)  _{n\times n}$,%

\[%
\begin{array}
[c]{c}%
\left(
\begin{array}
[c]{ccccc}%
1 & 0 & 0 & 0 & 0\\
1 & 1 & 0 & 0 & 0\\
2 & 2 & 1 & 0 & 0\\
4 & 5 & 3 & 1 & 0\\
9 & 12 & 9 & 4 & 1
\end{array}
\allowbreak\right)  ^{-1}=\left(
\begin{array}
[c]{ccccc}%
1 & 0 & 0 & 0 & 0\\
-1 & 1 & 0 & 0 & 0\\
0 & -2 & 1 & 0 & 0\\
1 & 1 & -3 & 1 & 0\\
-1 & 2 & 3 & -4 & 1
\end{array}
\allowbreak\right)  =\left(  m_{i,j}\right)  _{4\times4}\\
\text{Inverse Motzkin matrix when }\omega=1
\end{array}
\]
where%
\[
\sum_{i\geq0}m_{i,j}t^{i}=t^{j}\phi\left(  t\right)  ^{-j-1}%
\]
Note that $\left(  m_{i,j}\right)  $ is also a \emph{Riordan matrix}. The
above generating function for $m_{i,j}$ implies that%
\[
m_{i,j}=\left[  t^{i}\right]  \frac{1}{1+\omega t+t^{2}}\left(  \frac
{t}{1+\omega t+t^{2}}\right)  ^{j}=\left[  t^{i-j}\right]  \left(  1+\omega
t+t^{2}\right)  ^{-j-1}=C_{i-j}^{j+1}\left(  -\omega/2\right)  .
\]

The polynomials $C_{n}^{\lambda}\left(  x\right)  =\sum_{k=0}^{n/2}%
\binom{n-k+\lambda-1}{n-k}\binom{n-k}{n-2k}\left(  -1\right)  ^{k}\left(
2x\right)  ^{n-2k}$ are the \emph{Gegenbauer polynomials}, and therefore%
\begin{equation}
m_{i,j}=\sum_{l=0}^{\left(  i-j\right)  /2}\binom{i-l}{i-j-l}\binom{i-j-l}%
{l}\left(  -1\right)  ^{l}\left(  -\omega\right)  ^{i-j-2l}%
\label{(explictInverse)}%
\end{equation}
The recurrence relation for the (orthogonal) Gegenbauer polynomials%
\[
2x\left(  n+\lambda\right)  C_{n}^{\lambda}\left(  x\right)  =\left(
n+2\lambda-1\right)  C_{n-1}^{\lambda}+\left(  n+1\right)  C_{n+1}^{\lambda
}\left(  x\right)
\]
gives us immediately a recurrence for the inverse numbers $m_{i,j}$, $0\leq
j\leq i-1$,%
\[
\left(  i-j\right)  m_{i,j}=-\omega im_{i-1,j}-\left(  i+j\right)  m_{i-2,j}%
\]
with initial values $m_{i,j}=\delta_{i,j}$ for $j\geq i$.

We need later in the paper the following Motzkin ploynomial%
\begin{align}
\sum_{j=0}^{k}m_{k,j}t^{k-j}  & =\sum_{j=0}^{k}C_{j}^{k-j+1)}\left(
-\omega/2\right)  t^{j}\nonumber\\
& =\sum_{l=0}^{k/2}\sum_{j=0}^{k-2l}\binom{k-l}{k-j-l}\binom{k-j-l}%
{k-j-2l}\left(  -1\right)  ^{l}\left(  -\omega\right)  ^{k-j-2l}%
t^{k-j}\nonumber\\
& =\sum_{l=0}^{k/2}\binom{k-l}{l}\left(  -1\right)  ^{l}t^{2l}\left(  1-\omega
t\right)  ^{k-2l}\label{(MinversePolynomial)}%
\end{align}
$\allowbreak$

From%
\[
\left(  \left(  M\left(  i,j\right)  \right)  _{0\leq i,j\leq n}\right)
^{-1}=\left(  m_{i,j}\right)  _{0\leq i,j\leq n}%
\]
follows
\[
\sum_{k=0}^{n}M\left(  k,i;w\right)  m_{k,j}=\delta_{i,j}.
\]
However, in the case of Motzkin matrices more than this simple linear algebra
result holds.

\begin{lemma}
For all nonnegative integers $i$ and $i$ holds%
\[
M\left(  i,j;\omega\right)  =\sum_{k=0}^{j}m_{j,k}M_{i+k;\omega}%
\]
and%
\[
m_{i,j}=\sum_{k=0}^{i-j}m_{i+1,j+1+k}M_{k;\omega}%
\]

\end{lemma}

The proof can be done via generating functions. Note that
\[
\sum_{n\geq0}\sum_{j\geq0}x^{j}t^{n}M\left(  n,j;\omega\right)  =\frac
{\mu\left(  t\right)  }{1-xt\mu\left(  t\right)  }=\frac{1}{1+\omega
x+x^{2}-x/t}\left(  \mu\left(  t\right)  -\frac{x}{t}\right)
\]
and%
\[
\sum_{j\geq0}x^{j}\sum_{i\geq j}m_{i,j}t^{i}=\sum_{j\geq0}x^{j}t^{j}%
\phi\left(  t\right)  ^{j+1}=\frac{\phi\left(  t\right)  }{1-xt\phi\left(
t\right)  }=\frac{1}{1/\phi\left(  t\right)  -xt}=\frac{1}{1+\omega
t+t^{2}-xt}.
\]
Replace $t$ by $x$ and $x$ by $1/t$ in the above generating function for the
inverse $m_{i,j}$ to get the Laurent series%
\[
\sum_{j\geq0}t^{-j}\sum_{i\geq j}m_{i,j}x^{i}=\frac{1}{1+\omega x+x^{2}-x/t}%
\]
hence%
\[
\sum_{n\geq0}\sum_{j\geq0}x^{j}t^{n}M\left(  n,j;\omega\right)  =\left(
\mu\left(  t\right)  -\frac{x}{t}\right)  \sum_{j\geq0}t^{-j}\sum_{i\geq
j}m_{i,j}x^{i}%
\]
Now both sides must be power series in $x$ and $t$. This condition gives the
Lemma. The Lemma also has the

\begin{corollary}
\label{CorOrth}%
\begin{equation}
\sum_{k=0}^{j}m_{j,k}M_{i+k,w}=\delta_{i,j}\text{ for }0\leq i\leq
j\label{(MotzkinGeneralizedOrth)}%
\end{equation}

\end{corollary}

because $M\left(  i,j;\omega\right)  =\delta_{i,j}$ for all $0\leq i\leq j$.

\section{Two applications of the inverse Motzkin matrix}

The Lemma says that%
\[
\left(  m_{i,j}\right)  _{0\leq i,j\leq n}\left(  M_{i+j;\omega}\right)
_{0\leq i,j\leq n}=\left(  M\left(  i,j;\omega\right)  \right)  _{0\leq
i,j\leq n}%
\]
which gives a direct way of calculcating the first \emph{Hankel determinant}
\begin{equation}
\det\left(  M_{i+j;\omega}\right)  _{0\leq i,j\leq n}=\frac{1}{\det\left(
m_{i,j}\right)  }\det\left(  M\left(  i,j;w\right)  \right)
=1\label{(Hankel1Motzkin)}%
\end{equation}
However, subsequent Hankel determinants are more complicated; we want to show
a way how to calculate a determinant proposed by Cameron and Yip
\cite{Cameron}. For a broader theory of Hankel determinants in lattice path
enumeration see \cite{Cigler}.

\subsection{The Hankel determinant $\left\vert \alpha M_{i+j;\omega}+\beta
M_{i+j+1;\omega}\right\vert _{0\leq i,j\leq n-1}$}

The Hankel determinant of $\left(  \alpha M_{i+j;\omega}+\beta M_{i+j+1;\omega
}\right)  _{0\leq i,j\leq n-1}$ equals for $\omega=1$%
\begin{align*}
& =\left\vert
\begin{array}
[c]{ccccc}%
\alpha+2\beta & 2\alpha+4\beta & 4\alpha+7\beta & \dots & \alpha
M_{n-1;1}+\beta M_{n;1}\\
2\alpha+4\beta & 4\alpha+7\beta & 7\alpha+9\beta &  & \\
4\alpha+7\beta & 7\alpha+9\beta & 4\alpha+7\beta & \vdots & \vdots\\
7\alpha+9\beta & 9\alpha+21\beta & 9\alpha+21\beta &  & \\
\vdots & \vdots & \vdots &  & \\
\alpha M_{n-1;1}+\beta M_{n;1} & \alpha M_{n;1}+\beta M_{n+1;1} & \alpha
M_{n+1;1}+\beta M_{n+2;1} &  & \alpha M_{2n-2;1}+\beta M_{2n,n;1}%
\end{array}
\right\vert \\
& =\left\vert \left(  M_{i+j;1}\right)  _{0\leq i,j\leq n-1}\right\vert
\left\vert
\begin{array}
[c]{ccccc}%
\alpha & 0 & 0 & \dots & -\beta m_{n}\left(  0\right) \\
\beta & \alpha & 0 &  & -\beta m_{n}\left(  1\right) \\
0 & \beta & \alpha &  & -\beta m_{n}\left(  2\right) \\
&  &  & \vdots & \\
0 & 0 & 0 & \alpha & -\beta m_{n}\left(  n-2\right) \\
0 & 0 & 0 & \beta & \alpha-\beta m_{n}\left(  n-1\right)
\end{array}
\right\vert
\end{align*}
\newline because the last column in the matrix on the right when multiplied
with the $i$-th row of the matrix on the left gives $\alpha M_{i+n-1;\omega
}-\beta\sum_{k=0}^{n-1}m_{n,k}M_{i+k;\omega}=\alpha M_{i+n-1;\omega}+\beta
M_{i+n;\omega}-\beta\delta_{i,n}$ by Corollary \ref{CorOrth}. Now%

\begin{align*}
& \left\vert
\begin{array}
[c]{ccccc}%
\alpha & 0 & 0 & \dots & -\beta m_{n,0}\\
\beta & \alpha & 0 &  & -\beta m_{n,1}\\
0 & \beta & \alpha & \dots & -\beta m_{n,2}\\
&  &  & \vdots & \\
0 & 0 & \dots & \alpha & -\beta m_{n,n-2}\\
0 & 0 & \dots & \beta & \alpha-\beta m_{n,n-1}%
\end{array}
\right\vert \\
& =\alpha^{-\binom{n}{2}}\left\vert
\begin{array}
[c]{ccccc}%
\alpha & 0 & 0 & \dots & -\beta m_{n,0}\\
\alpha\beta & \alpha^{2} & 0 &  & -\alpha\beta m_{n,1}\\
0 & \alpha^{2}\beta & \alpha^{3} & \dots & -\alpha^{2}\beta m_{n,2}\\
&  &  & \vdots & \\
0 & 0 & \dots & \alpha^{n-1} & -\alpha^{n-2}\beta m_{n,n-2}\\
0 & 0 & \dots & \alpha^{n-1}\beta & \alpha^{n-1}-\alpha^{n-1}\beta m_{n,n-1}%
\end{array}
\right\vert \\
& =\alpha^{-\binom{n}{2}}\left\vert
\begin{array}
[c]{ccccc}%
\alpha & 0 & 0 & \dots & -\beta m_{n,0}\\
0 & \alpha^{2} & 0 &  & \beta^{2}m_{n,0}-\alpha\beta m_{n,1}\\
0 & 0 & \alpha^{3} & \dots & -\beta^{3}m_{n,0}+\alpha\beta^{2}m_{n,1}%
-\alpha^{2}\beta m_{n,2}\\
&  &  & \vdots & \\
0 & 0 & \dots & \alpha^{n-1} & -\sum_{i=0}^{n-2}\left(  -1\right)
^{n-2-i}\beta^{n-1-i}\alpha^{i}m_{n,i}\\
0 & 0 & \dots & 0 & \alpha^{n}-\sum_{i=0}^{n-1}\left(  -1\right)
^{n-1-i}\beta^{n-i}\alpha^{i}m_{n,i}%
\end{array}
\right\vert
\end{align*}
Therefore $\det\left(  \left(  \alpha M_{i+j;\omega}+\beta M_{i+j+1;\omega
}\right)  _{0\leq i,j\leq n-1}\right)  =\alpha^{n}-\sum_{i=0}^{n-1}\left(
-1\right)  ^{n-1-i}\beta^{n-i}\alpha^{i}m_{n,i}=\sum_{i=0}^{n}\left(
-\beta\right)  ^{n-i}\alpha^{i}m_{n,i}=\sum_{i=0}^{n}\left(  -1\right)
^{n-i}\beta^{n-i}\alpha^{i}P_{n-i}^{(-i-1)}\left(  -\omega/2\right)  $. This
can be written explicitly as $\det\left(  \left(  \alpha M_{i+j;\omega}+\beta
M_{i+j+1;\omega}\right)  _{0\leq i,j\leq n-1}\right)  =$%

\begin{align*}
& \left(  -\beta\right)  ^{n}\sum_{k=0}^{n}\left(  -\alpha/\beta\right)
^{k}m_{n,k}\\
& =\left(  -\beta\right)  ^{n}U_{n}\left(  \frac{-\alpha/\beta-\omega}%
{2}\right)  =\left(  -\beta\right)  ^{n}\sum_{k=0}^{n/2}\binom{n-k}{k}\left(
-1\right)  ^{k}\left(  -\alpha/\beta-\omega\right)  ^{n-2k}\\
& =\sum_{k=0}^{n/2}\binom{n-k}{k}\left(  -1\right)  ^{k}\beta^{2k}\left(
\alpha+\beta\omega\right)  ^{n-2k}\\
& =\frac{2^{-n-1}}{\sqrt{\left(  \alpha+\omega\beta\right)  ^{2}-4\beta^{2}}%
}\times\\
& \times\left(  \left(  \sqrt{\left(  \alpha+\omega\beta\right)  ^{2}%
-4\beta^{2}}+\alpha+\omega\beta\right)  ^{n+1}+\left(  \sqrt{\left(
\alpha+\beta\omega\right)  ^{2}-4\beta^{2}}-\alpha-\beta\omega\right)
^{n+1}\right)
\end{align*}
If $\alpha=\beta=1$, then $\det\left(  \left(  M_{i+j;\omega}+M_{i+j+1;\omega
}\right)  _{0\leq i,j\leq n-1}\right)  =$%
\begin{align*}
& \frac{1}{2^{n+1}\sqrt{\left(  \omega+1\right)  ^{2}-4}}\left(  \left(
1+\omega+\sqrt{\left(  \omega+1\right)  ^{2}-4}\right)  ^{n+1}-\left(
1+\omega-\sqrt{\left(  \omega+1\right)  ^{2}-4}\right)  ^{n+1}\right) \\
& =\sum_{k=0}^{n}\left(  -1\right)  ^{n-k}\binom{k}{n-k}\left(  \omega
+1\right)  ^{2k-n}%
\end{align*}

which approaches $n+1$ if $\omega\rightarrow1$. In the case of Dyck path, we
obtain $\delta_{0,n}$ for this determinat of the sum of matrices. If $\beta=1$
and $\alpha=0$, then the determinant is the second Hankel determinant of the
Motzkin numbers,%
\[
\det\left(  \left(  M_{i+j+1;\omega}\right)  _{0\leq i,j\leq n-1}\right)
=\sum_{k=0}^{n/2}\binom{n-k}{k}\left(  -1\right)  ^{k}\omega^{n-2k}%
\]
If $\alpha=1$ and $\beta=0$ then $\det\left(  M_{i+j;\omega}\right)  _{0\leq
i,j\leq n-1}=1$, independent of $\omega$ (see (\ref{(Hankel1Motzkin)}). The
same approach also shows the recursion
\[
\left\vert M_{i+j+2;\omega}\right\vert _{0\leq i,j\leq n-1}=\left\vert
M_{i+j+2;\omega}\right\vert _{0\leq i,j\leq n-2}+\left\vert M_{i+j+1;\omega
}\right\vert _{0\leq i,j\leq n-1}^{2}%
\]

\subsection{Motzkin in a band}

The number of Motzkin paths staying strictly below the line $y=k$ for $k>0$ is
known to have the generating function \cite[Proposition 12]{Flajolet}%
\[
\sum_{n\geq0}M_{n}^{\left(  k\right)  }t^{n}=\mu\left(  t\right)
\frac{1-\left(  t\mu\left(  t\right)  \right)  ^{2k}}{1-\left(  t\mu\left(
t\right)  \right)  ^{2\left(  k+1\right)  }}=\frac{1}{t}\frac{\left(  \frac
{1}{t\mu}\right)  ^{k}-\left(  t\mu\right)  ^{k}}{\left(  \frac{1}{t\mu
}\right)  ^{k+1}-\left(  t\mu\right)  ^{k+1}}%
\]

\[%
\begin{tabular}
[c]{l||llllllllll}%
$k$ & 0 & 0 & 0 & 0 & 0 & 0 & 0 & 0 & 0 & 0\\
3 &  &  &  & $1$ & $4$ & $14$ & $44$ & $133$ & $392$ & $1140$\\
2 &  &  & $1$ & $3$ & $9$ & $25$ & $69$ & $189$ & $518$ & $1422$\\
1 &  & $1$ & $2$ & $5$ & $12$ & $30$ & $76$ & $196$ & $512$ & $1353$\\
$\boldsymbol{0}$ & $1$ & $1$ & $2$ & $4$ & $9$ & $21$ & $51$ & $127$ & $323$ &
$835$\\
-1 & 0 & 0 & 0 & 0 & 0 & 0 & 0 & 0 & 0 & 0\\\hline\hline
& 0 & 1 & 2 & 3 & 4 & 5 & 6 & 7 & 8 & $n$\\
\multicolumn{11}{c}{$M_{n}^{(4)}$ is given in row $0.$}%
\end{tabular}
\]
From $\mu\left(  t\right)  \left(  1-\omega t\right)  -1=t^{2}\mu\left(
t\right)  ^{2}$ (see (\ref{(quadratic)})) follows%
\begin{align*}
\mu_{1,2}\left(  t\right)   & =\left(  1-\omega t\pm\sqrt{\left(  1-\omega
t\right)  ^{2}-4t^{2}}\right)  /\left(  2t^{2}\right) \\
t\mu_{1,2}\left(  t\right)   & =\left(  1-\omega t\pm\sqrt{\left(  1-\omega
t\right)  ^{2}-4t^{2}}\right)  /\left(  2t\right)
\end{align*}
thus%
\[
\mu_{1}+\mu_{2}=\left(  1-\omega t\right)  /t^{2}\text{ and }\mu_{1}\mu
_{2}=1/t^{2}%
\]

Hence
\begin{align*}
\sum_{n\geq0}M_{n;\omega}^{\left(  k\right)  }t^{n}  & =\frac{1}{t}%
\frac{\left(  t\mu_{1}\right)  ^{k}-\left(  t\mu_{2}\right)  ^{k}}{\left(
t\mu_{1}\right)  ^{k+1}-\left(  t\mu_{2}\right)  ^{k+1}}=\frac{1}{t}%
\frac{\left(  t\mu_{2}\right)  ^{-k}-\left(  t\mu_{1}\right)  ^{-k}}{\left(
t\mu_{2}\right)  ^{-k-1}-\left(  t\mu_{1}\right)  ^{-k-1}}\\
& =\frac{\sum_{j=0}^{\left(  k-1\right)  /2}\left(  -1\right)  ^{j}%
\binom{k-1-j}{j}t^{2j}\left(  1-\omega t\right)  ^{k-1-2j}}{\sum_{j=0}%
^{k/2}\left(  -1\right)  ^{j}\binom{k-j}{j}t^{2j}\left(  1-\omega t\right)
^{k-2j}}\\
& =\frac{\sum_{i=0}^{k-1}m_{k-1,i}t^{k-1-i}}{\sum_{i=0}^{k}m_{k,i}t^{k-i}}%
\end{align*}
(see (\ref{(MinversePolynomial)})). The OEIS lists many special cases for $k$;
here are a few, with $\omega=1$.

\begin{enumerate}
\item $\sum_{n\geq0}M_{n;1}^{\left(  1\right)  }t^{n}\frac{1}{1-t}%
\iff1,1,1,1,\dots$

\item $\sum_{n\geq0}M_{n;1}^{\left(  2\right)  }t^{n}=\frac{1-t}{\left(
1-t\right)  ^{2}-t^{2}}=1+t+2t^{2}+4t^{3}+8t^{4}+16t^{5}\dots$\newline thus
$1,1,2,4,8,16,32,64,\dots$, the powers of $2$.

\item $\sum_{n\geq0}M_{n;1}^{\left(  3\right)  }t^{n}=\frac{2t-1}{\left(
1-t\right)  \left(  t^{2}+2t-1\right)  }\ $thus $1,1,2,4,9,21,50,120,\dots$ (A171842)

\item $\sum_{n\geq0}M_{n;1}^{\left(  4\right)  }t^{n}=\left(  1-3t+t^{2}%
+t^{3}\right)  /\left(  1-4t+3t^{2}+2t^{3}-t^{4}\right)  $, thus\newline%
$1,1,2,4,9,21,51,127,322,826,\dots:$ (A005207), generating function by Alois
P. Heinz.
\end{enumerate}

The special form of the generating function%
\begin{equation}
\sum_{n\geq0}M_{n;\omega}^{\left(  k\right)  }t^{n}=\frac{\sum_{i=0}%
^{k-1}m_{k-1,i}t^{k-1-i}}{\sum_{i=0}^{k}m_{k,i}t^{k-i}}\label{(BoundedM)}%
\end{equation}
works with weight $\omega$, for all $k=1,2,\dots$.It is equivalent to the
recursion $\sum_{j=0}^{k}M_{n-j}^{\left(  k\right)  }m_{k,k-j}=0$ for all
$n\geq k$, with initial values $\sum_{j=0}^{n}M_{n-j}^{\left(  k\right)
}m_{k,k-j}=m_{k-1,k-1-n}$ for all $n=0,\dots,k-1$.

\section{Horizontal steps of length $w$}

A \textquotedblleft natural\textquotedblright\ generalization of Motzkin paths
is a lattice path $W$ that takes horizontal steps of some positive length $w$,
weighted by $\omega$. We would like to see similar results as
(\ref{(BoundedM)}) in such cases. However, we have a result only for the case
$w=2$, the Schr\"{o}der paths.%

\[%
\begin{tabular}
[c]{l||llllllllll}%
$\uparrow m$ &  &  &  &  &  &  &  &  & $1$ & $0$\\
7 &  &  &  &  &  &  &  & $1$ & $0$ & $8$\\
6 &  &  &  &  &  &  & $1$ & $0$ & $7$ & $7\omega$\\
5 &  &  &  &  &  & $1$ & $0$ & $6$ & $6\omega$ & $27$\\
4 &  &  &  &  & $1$ & $0$ & $5$ & $5\omega$ & $20$ & $35\omega$\\
3 &  &  &  & $1$ & $0$ & $4$ & $4\omega$ & $14$ & $24\omega$ & $48+10\omega
^{2}$\\
2 &  &  & $1$ & $0$ & $3$ & $3\omega$ & $9$ & $15\omega$ & $28+6\omega^{2} $ &
$63\omega$\\
1 &  & $1$ & $0$ & $2$ & $2\omega$ & $5$ & $8\omega$ & $14+3\omega^{2}$ &
$30\omega$ & $42+20\omega^{2}$\\
0 & $1$ & $0$ & $1$ & $\omega$ & $2$ & $3\omega$ & $5+\omega^{2}$ & $10\omega$
& $14+6\omega^{2}$ & $35\omega+\omega^{3}$\\\hline\hline
& 0 & 1 & 2 & 3 & 4 & 5 & 6 & 7 & 8 & $n\rightarrow$\\
\multicolumn{11}{c}{$w=3$ ($W_{n}$ is given in row $0$)}%
\end{tabular}
\]

\subsection{The recursion for $W$}

Let us consider the step set $\left\{  \nearrow,\searrow,\longrightarrow
^{w}\right\}  $, where $\rightarrow^{w}=:\left(  w,0\right)  $, for any
positive integer $w$. Denote the number of paths from $\left(  0,0\right)  $
to $\left(  n,j\right)  $ by $W\left(  n,j;\omega\right)  $, where the
horizontal steps are weighted by $\omega$. We get the recursion
\begin{align*}
W\left(  n,j;\omega\right)   & =W\left(  n-1,j+1;\omega\right)  +W\left(
n-1,j-1;\omega\right)  +\omega W\left(  n-w,j;\omega\right) \\
W\left(  n,j;\omega\right)   & =0\text{ for }j<n\\
W_{n;\omega}  & =W(n,0;\omega).
\end{align*}
The generating function is well known,
\begin{equation}
\sum_{n\geq0}W_{n;\omega}t^{n}=\frac{1-\omega t^{w}-\sqrt{\left(  1-\omega
t^{w}\right)  ^{2}-4t^{2}}}{2t^{2}}=:\mu_{w}\left(  t;\omega\right)
\label{(W(t,0))}%
\end{equation}

The recursion can be reformulated as%
\[
W\left(  n,j;\omega\right)  =W\left(  n+1,j-1;\omega\right)  -W\left(
n,j-2;\omega\right)  -\omega W\left(  n+1-w,j-1;\omega\right)  \text{ for
}m\geq n
\]
We find the generating function identity $\sum_{i\geq0}W\left(  i,j;\omega
\right)  t^{i}=$%
\begin{align*}
& \sum_{i\geq0}W\left(  i+1,j-1;\omega\right)  t^{i}-\sum_{i\geq0}W\left(
i,j-2;\omega\right)  t^{i}\\
& -\omega\sum_{i\geq w-1}W\left(  i+1-w,j-1;\omega\right)  t^{i}\\
& =\sum_{i\geq0}W\left(  i+1,j-1;\omega\right)  t^{i}-\omega\left(
\sum_{i\geq-1}W\left(  i+1,j-1;\omega\right)  t^{i+1+w-1}\right) \\
& -\sum_{i\geq0}W\left(  i,j-2;\omega\right)  t^{i}\\
& =t^{-1}\sum_{i\geq0}W\left(  i+1,j-1;\omega\right)  t^{i+1}-\omega
t^{w-1}\left(  \sum_{i\geq-1}W\left(  i+1,j-1;\omega\right)  t^{i+1}\right) \\
& -\sum_{i\geq0}W\left(  i,j-2;\omega\right)  t^{i}\\
& =\left(  t^{-1}-\omega t^{w-1}\right)  \left(  \sum_{i\geq0}W\left(
i,j-1;\omega\right)  t^{i}-\delta_{j,1}\right)  -\sum_{i\geq0}W\left(
i,j-2;\omega\right)  t^{i}%
\end{align*}
Let $\mathcal{W}\left(  t,j;\omega\right)  =\sum_{i\geq0}W\left(
i,j;\omega\right)  t^{i}$. In this notation,%
\begin{align}
\mathcal{W}\left(  t,j;\omega\right)   & =\frac{1-\omega t^{w}}{t}%
\mathcal{W}\left(  t,j-1;\omega\right)  -\mathcal{W}\left(  t,j-2;\omega
\right)  \text{ for }j>1\label{(Wrecursion)}\\
\mathcal{W}\left(  t,1;\omega\right)   & =\frac{1}{t}\left(  \left(  1-\omega
t^{w}\right)  \mathcal{W}\left(  t,0;\omega\right)  -1\right) \nonumber
\end{align}
For example,

$\mathcal{W}\left(  t,2;\omega\right)  =\frac{\left(  1-\omega t^{w}\right)
}{t}\mathcal{W}\left(  t,1;\omega\right)  -\mathcal{W}\left(  t,0;\omega
\right)  =\frac{\left(  1-\omega t^{w}\right)  }{t}\frac{1}{t}\left(  \left(
1-\omega t^{w}\right)  \mathcal{W}\left(  t,0;\omega\right)  -1\right)
-\mathcal{W}\left(  t,0;\omega\right)  $\newline$=\left(  \frac{\left(
1-\omega t^{w}\right)  ^{2}}{t^{2}}-1\right)  \mathcal{W}\left(
t,0;\omega\right)  -\frac{\left(  1-\omega t^{w}\right)  }{t^{2}}$, and
$\mathcal{W}\left(  t,0;\omega\right)  =\mu_{w}\left(  t;\omega\right)  $ is
given in (\ref{(W(t,0))}).

$\mathcal{W}\left(  t,3;\omega\right)  =\frac{\left(  1-\omega t^{w}\right)
}{t}\mathcal{W}\left(  t,2;\omega\right)  -\mathcal{W}\left(  t,1;\omega
\right)  $\newline$=\frac{\left(  1-\omega t^{w}\right)  }{t}\left(  \left(
\frac{\left(  1-\omega t^{w}\right)  ^{2}}{t^{2}}-1\right)  \mathcal{W}\left(
t,0;\omega\right)  -\frac{\left(  1-\omega t^{w}\right)  }{t^{2}}\right)
-\frac{1}{t}\left(  \left(  1-\omega t^{w}\right)  \mathcal{W}\left(
t,0;\omega\right)  -1\right)  $\newline$=\left(  \frac{\left(  1-\omega
t^{w}\right)  ^{2}}{t^{2}}-2\right)  \left(  \frac{\left(  1-\omega
t^{w}\right)  }{t}\mu_{w}\left(  t;\omega\right)  \right)  +\frac{1}{t}%
-\frac{\left(  1-\omega t^{w}\right)  ^{2}}{t^{3}}$

We find an explicit expression for $\mathcal{W}\left(  t,j;\omega\right)  $ in
the next section.

\subsection{Solution to Recursion for $\mathcal{W}$ and $\mathcal{W}^{\left(
k\right)  }$}

The linear recursion (\ref{(Wrecursion)}) is called Fibonacci-like. It is of
the form
\[
\sigma_{n}=u\sigma_{n-1}+v\sigma_{n-2}%
\]
with $u=\frac{1-\omega t^{w}}{t}$ and $v=-1$, for $n>1$. We know the inital
values $\sigma_{0}$ and $\sigma_{1}=u\sigma_{0}-1/t$.

Hence $\sigma_{n}=\left[  \tau^{n}\right]  \frac{\sigma_{0}+\left(  \sigma
_{1}-u\sigma_{0}\right)  \tau}{1-u\tau-v\tau^{2}}=\left[  \tau^{n}\right]
\frac{\sigma_{0}-\tau/t}{1-u\tau+\tau^{2}}$ in this case, or $\sigma_{n}=$%
\begin{align}
& \left[  \tau^{n}\right]  \left(  \sigma_{0}-\frac{\tau}{t}\right)
\sum_{i=0}^{\infty}\binom{i}{j}\left(  -1\right)  ^{j}u^{i-j}\tau
^{i-j+2j}\label{(Solution to Rec)}\\
& =\sigma_{0}\sum_{j=0}^{n}\tbinom{n-j}{j}\left(  -1\right)  ^{j}\left(
\frac{1-\omega t^{w}}{t}\right)  ^{n-2j}-\frac{1}{t}\sum_{j=0}^{n-1}%
\tbinom{n-1-j}{j}\left(  -1\right)  ^{j}\left(  \frac{1-\omega t^{w}}%
{t}\right)  ^{n-1-2j}\nonumber
\end{align}

Let us define
\begin{equation}
p_{n}\left(  t\right)  :=\sum_{j=0}^{n}\tbinom{n-j}{j}\left(  \frac{1-\omega
t^{w}}{t}\right)  ^{n-2j}\left(  -1\right)  ^{j}\label{(p_n)}%
\end{equation}
Hence%
\begin{equation}
\mathcal{W}\left(  t,j;\omega\right)  =\left(  \frac{1-\omega t^{w}%
-\sqrt{\left(  1-\omega t^{w}\right)  ^{2}-4t^{2}}}{2t^{2}}\right)
p_{j}\left(  t\right)  -p_{j-1}\left(  t\right)  /t\label{(W(t,j))}%
\end{equation}
where $p_{j}=0$ for all $j<0$.

The generating function $\mathcal{W}^{\left(  k\right)  }\left(
t,j;\omega\right)  =\sum_{n\geq0}^{\left(  k\right)  }W^{\left(  k\right)
}\left(  n,j;\omega\right)  t^{n}$ is generating the case where the lattice
paths stay strictly below $y=k$; the numbers $W^{\left(  k\right)  }\left(
n,j;\omega\right)  $ are the number of paths with $\omega$-weighted horizontal
steps of length $w$, and diagonal up and down steps, that do not reach the
line $y=k$, and stay above the $x$-axis. That means, $0\leq j<k$. We also know
$\mathcal{W}^{\left(  k\right)  }\left(  t,0;\omega\right)  $
\begin{align*}
& =\sum_{n\geq0}W_{n}^{\left(  k\right)  }t^{n}=\mu_{w}\left(  t;\omega
\right)  \frac{1-\left(  t\mu_{w}\left(  t;\omega\right)  \right)  ^{2k}%
}{1-\left(  t\mu_{w}\left(  t;\omega\right)  \right)  ^{2\left(  k+1\right)
}}\\
& =\frac{1-\omega t^{w}-\sqrt{\left(  1-\omega t^{w}\right)  ^{2}-4t^{2}}%
}{2t^{2}}\frac{1-\left(  \frac{1-\omega t^{w}-\sqrt{\left(  1-\omega
t^{w}\right)  ^{2}-4t^{2}}}{2t}\right)  ^{2k}}{1-\left(  \frac{1-\omega
t^{w}-\sqrt{\left(  1-\omega t^{w}\right)  ^{2}-4t^{2}}}{2t}\right)
^{2\left(  k+1\right)  }}%
\end{align*}

The recursion is the same as for $\mathcal{W}\left(  t,j;\omega\right)  $.
Only the initial values have changed (see $\mathcal{W}^{\left(  k\right)
}\left(  t,0;\omega\right)  $ above).

We get%
\begin{align*}
& \mathcal{W}^{\left(  k\right)  }\left(  t,j;\omega\right) \\
& =\left(  \mu_{w}\left(  t;\omega\right)  \frac{1-\left(  t\mu_{w}\left(
t;\omega\right)  \right)  ^{2k}}{1-\left(  t\mu_{w}\left(  t;\omega\right)
\right)  ^{2\left(  k+1\right)  }}\right)  p_{j}\left(  t\right)
-p_{j-1}\left(  t\right)
\end{align*}
and $\sum_{n\geq0}W_{n}^{\left(  k\right)  }t^{n}$%

\begin{align}
& =\frac{1}{2t^{2}}\frac{\left(  \frac{1-\omega t^{w}-\sqrt{\left(  1-\omega
t^{w}\right)  ^{2}-4t^{2}}}{2t}\right)  ^{-k}-\left(  \frac{1-\omega
t^{w}-\sqrt{\left(  1-\omega t^{w}\right)  ^{2}-4t^{2}}}{2t}\right)  ^{k}%
}{\left(  \frac{1-\omega t^{w}-\sqrt{\left(  1-\omega t^{w}\right)
^{2}-4t^{2}}}{2t}\right)  ^{-\left(  k+1\right)  }-\left(  \frac{1-\omega
t^{w}-\sqrt{\left(  1-\omega t^{w}\right)  ^{2}-4t^{2}}}{2t}\right)  ^{\left(
k+1\right)  }}\nonumber\\
& =2\frac{\sum_{i=0}^{\left(  k-1\right)  /2}\left(  -1\right)  ^{i}%
2^{2i}t^{2i}\left(  1-\omega t^{w}\right)  ^{k-2i}\sum_{j=0}^{\infty}%
\binom{j+i}{i}\binom{k}{2j+2i+1}}{\sum_{i=0}^{\left(  k+1\right)  /2}\left(
-1\right)  ^{i}2^{2i}t^{2i}\left(  1-\omega t^{w}\right)  ^{k+1-2i}\sum
_{j=0}^{\infty}\binom{j+i}{i}\binom{k+1}{2j+1+2i}}\nonumber\\
& =\frac{\sum_{i=0}^{k/2}\left(  -1\right)  ^{i}t^{2i}\left(  1-\omega
t^{w}\right)  ^{k-2i}\binom{k-i-1}{i}}{\sum_{i=0}^{\left(  k+1\right)
/2}\left(  -1\right)  ^{i}t^{2i}\left(  1-\omega t^{w}\right)  ^{k+1-2i}%
\binom{k-i}{i}}=\frac{p_{k-1}\left(  t\right)  }{tp_{k}\left(  t\right)
}\label{(genfunW_n^(k))}%
\end{align}
where $p_{n}\left(  t\right)  $ is given in (\ref{(p_n)}).

\section{Schr\"{o}der numbers}

If $w=2$, then every horizontal steps gains two units. We denote the number of
paths to $\left(  n,j\right)  $ by $S\left(  n,j;\omega\right)  $, the general
weighted Schr\"{o}der numbers.
\[%
\begin{array}
[c]{c}%
\begin{array}
[c]{cccccccc}%
j\uparrow &  &  &  &  & 1 & 0 & 5+5\omega\\
3 &  &  &  & 1 & 0 & 4+4\omega & 0\\
2 &  &  & 1 & 0 & 3+3\omega & 0 & 9+15\omega+6\omega^{2}\\
1 &  & 1 & 0 & 2+2\omega & 0 & 5+8\omega+3\omega^{2} & 0\\
0 & 1 & 0 & 1+\omega & 0 & 2+3\omega+\omega^{2} & 0 & 5+10\omega+6\omega
^{2}+\omega^{3}\\
n\rightarrow & 0 & 1 & 2 & 3 & 4 & 5 & 6
\end{array}
\\
\text{The general weighted Schr\"{o}der numbers }S\left(  n,j;\omega\right)
\text{. The numbers }S_{n;\omega}\text{ are in row }0\text{.}%
\end{array}
\]
The following matrix contains the \textquotedblleft\emph{compressed}%
\textquotedblright\ Schr\"{o}der numbers by removing the zeroes and shifting
all entries into the empty places. This is the same effect as replacing
$t^{2}$ in $\sum_{n=0}^{\infty}S\left(  n,j;\omega\right)  t^{j}$ by $t$.%
\[%
\begin{array}
[c]{c}%
\left(
\begin{array}
[c]{ccccc}%
1 & 0 & 0 & 0 & 0\\
2 & 1 & 0 & 0 & 0\\
6 & 4 & 1 & 0 & 0\\
22 & 16 & 6 & 1 & 0\\
90 & 68 & 30 & 8 & 1
\end{array}
\right)  ^{-1}=\left(
\begin{array}
[c]{ccccc}%
1 & 0 & 0 & 0 & 0\\
-2 & 1 & 0 & 0 & 0\\
2 & -4 & 1 & 0 & 0\\
-2 & 8 & -6 & 1 & 0\\
2 & -12 & 18 & -8 & 1
\end{array}
\right) \\
\text{compressed Schr\"{o}der numbers (}\omega=1\text{)~~~~ inverse compressed
Schr\"{o}der numbers}%
\end{array}
\]
The power series $\sum_{n=0}^{\infty}S\left(  n,j;\omega\right)  t^{j}$ is
given in (\ref{(W(t,0))}). For the compressed Schr\"{o}der numbers this
equation says%
\[
\mathcal{S}\left(  t,j;\omega\right)  =\sum_{n=0}^{\infty}S\left(
n,j;\omega\right)  t^{j}=\left(  \frac{1-\omega t-\sqrt{\left(  1-\omega
t\right)  ^{2}-4t}}{2t}\right)  p_{j}\left(  t\right)  -p_{j-1}\left(
t\right)  /t
\]
where
\begin{equation}
p_{n}\left(  t\right)  =t^{-n}\sum_{j=0}^{n}\tbinom{n-j}{j}t^{j}\left(
1-\omega t\right)  ^{n-2j}\left(  -1\right)  ^{j}\label{(cp-p_n)}%
\end{equation}

All references to Schr\"{o}der numbers will from now on mean the compressed
Schr\"{o}der numbers. Note that%

\begin{equation}
\mathcal{S}^{\left(  k\right)  }\left(  t;\omega\right)  =\sum_{n\geq0}%
S_{n}^{\left(  k\right)  }t^{n}=\frac{p_{k-1}\left(  t\right)  }{tp_{k}\left(
t\right)  }\label{(GenFunS_n^(k))}%
\end{equation}
by (\ref{(genfunW_n^(k))}).

\subsection{Inverse Schr\"{o}der numbers}

From (\ref{(W(t,0))}) we see that%

\[
\mu_{s}\left(  t\right)  =1+\omega t^{2}\mu_{s}\left(  t\right)  +t^{2}\mu
_{s}\left(  t\right)  ^{2}.
\]
Hence
\begin{align*}
\phi\left(  t\mu_{s}\left(  t\right)  \right)    & =\mu_{s}\left(  t\right)
=1+\omega t^{2}\mu_{s}\left(  t\right)  +t^{2}\mu_{s}\left(  t\right)  ^{2}\\
\phi\left(  t\right)    & =1+\frac{\omega t^{2}}{\phi\left(  t\right)  }+t^{2}%
\end{align*}
thus $\phi\left(  t\right)  =\frac{1}{2}+\frac{1}{2}t^{2}+\frac{1}{2}%
\sqrt{\left(  1+t^{2}\right)  ^{2}+4t^{2}\omega}$, a power series in $t^{2}$.
We let $\xi=t^{2}$ and get%
\begin{align*}
\phi\left(  \xi\right)    & =\frac{1}{2}+\frac{1}{2}\xi+\frac{1}{2}%
\sqrt{\left(  1+\xi\right)  ^{2}+4\xi\omega}\\
\mu_{s}\left(  \xi\right)    & =1+\omega\xi\mu_{s}\left(  t\right)  +\xi
\mu_{s}\left(  \xi\right)  ^{2}\\
& =\frac{1-\omega\xi-\sqrt{\left(  1-\omega\xi\right)  ^{2}-4\xi}}{2\xi}%
\end{align*}

Lagrange inversion tells us that for all $0\leq i\leq k$ holds%
\[
\left(  i+1\right)  \left[  \mu_{s}^{-k-1}\right]  _{k-i}=\left(  k+1\right)
\left[  \phi^{-i-1}\right]  _{k-i}=\left(  k+1\right)  s_{k,i}%
\]
$\allowbreak$and therefore%
\begin{align*}
s_{k,j}  & =\left[  \mu_{s}^{-k-1}\right]  _{k-j}=\frac{j+1}{k+1}\left[
t^{k-j}\right]  \left(  \frac{1-\omega t-\sqrt{\left(  1-\omega t\right)
^{2}-4t}}{2t}\right)  ^{-k-1}\\
& =\left[  t^{k-j}\right]  \frac{j+1}{k+1}\left(  \frac{2t\left(  1-\omega
t+\sqrt{\left(  1-\omega t\right)  ^{2}-4t}\right)  }{4t}\right)  ^{k+1}\\
& =\left[  t^{k-j}\right]  \frac{j+1}{k+1}\left(  \frac{1}{2}\left(  1-\omega
t\right)  \left(  1+\sqrt{1-\frac{4t}{\left(  1-\omega t\right)  ^{2}}%
}\right)  \right)  ^{k+1}\\
& =\frac{j+1}{k+1}\left(  -1\right)  ^{k-j}\sum_{m=0}^{k-j}2^{2m-k-1}%
\dbinom{k+1-2m}{k-j-m}\omega^{k-j-m}\sum_{l=0}^{k+1}\dbinom{k+1}{l}%
\dbinom{l/2}{m}%
\end{align*}
Now $\sum_{l=0}^{k+1}\binom{k+1}{l}\binom{\frac{1}{2}l}{m}=\frac{k+1}%
{k-2m+1}\binom{k-m}{m}2^{k+1-2m}\allowbreak$ for $0\leq m$ (see Gould
\cite[(3.163)]{Gould}, who attributes the formula to Carlitz). Therefore%
\[
s_{k,j}=\left(  -1\right)  ^{k-j}\sum_{m=0}^{k-j}\dbinom{k+1-2m}{k-j-m}%
\frac{j+1}{k-m+1}\binom{k-m+1}{m}\omega^{k-j-m},
\]
the compressed weighted \emph{inverse Schr\"{o}der numbers}. We need the
following polynomials: $\sum_{k\geq0}s_{n,k}t^{n-k}$%

\begin{align*}
& =\sum_{k=0}^{n}\sum_{m=0}^{n-k}\frac{k+1}{n-m+1}\binom{n-m+1}{m}%
\dbinom{n+1-2m}{n-k-m}t^{n-k}\left(  -1\right)  ^{n-k}\omega^{n-k-m}\\
& =t^{n}\sum_{m=0}^{n}\frac{1}{n-m+1}\binom{n-m+1}{m}\omega^{-m}\sum
_{k=0}^{\infty}\left(  k+1\right)  \dbinom{n+1-2m}{n-k-m}\left(
-\omega\right)  ^{n-k}t^{-k}\\
& =\sum_{m=0}^{n}\frac{1}{n-m+1}\binom{n-m+1}{m}\left(  -1\right)
^{m+1}\left(  1-\omega t\right)  ^{n-2m}t^{m}\left(  tm\omega+m-n-1\right)
\end{align*}
Hence%
\begin{align}
s_{n}\left(  t\right)   & =\sum_{k\geq0}s_{n,k}t^{n-k}%
\label{(inverseSchroeder)}\\
& =\sum_{m=0}^{n}\left(  \frac{t\omega m}{n-m+1}-1\right)  \binom{n-m+1}%
{m}\left(  -1\right)  ^{m+1}\left(  1-\omega t\right)  ^{n-2m}t^{m}\nonumber
\end{align}

\[%
\begin{tabular}
[c]{c}%
$%
\begin{array}
[c]{ccccc}%
1 & 0 & 0 & 0 & 0\\
-2 & 1 & 0 & 0 & 0\\
2 & -4 & 1 & 0 & 0\\
-2 & 8 & -6 & 1 & 0\\
2 & -12 & 18 & -8 & 1
\end{array}
$\\
\multicolumn{1}{l}{The compressed inverse $\left(  s_{n,k}\right)  $ for
$\omega=1$}%
\end{tabular}
\]

This matrix is A080246 in the OEIS. At the same reference we find the
generating function of the $k$-th column,%
\[
\sum_{n\geq k}s_{n,k}t^{n}=\left(  \frac{1-t}{1+t}\right)  ^{k}.
\]
Also,
\[
\sum_{n\geq0}\sum_{k=0}^{n}s_{n,k}t^{n-k}=\sum_{k=0}^{\infty}t^{-k}\left(
\frac{1-t}{1+t}\right)  ^{k}=\frac{t\left(  1+t\right)  }{2t+t^{2}-1}.
\]
Example: $\sum_{k=0}^{n}s_{4,k}t^{4-k}=s_{4}\left(  t\right)  =1-8t+18t^{2}%
-12t^{3}+2t^{4}$.

\subsection{Delannoy numbers}

The numbers $D\left(  n,k\right)  =\sum_{l=0}^{n}\binom{k}{l}\binom{n+k-l}%
{k}\omega^{l}$ are the Delannoy numbers; the numbers $D\left(  n,n+j\right)  $
are counting all weighted Grand Schr\"{o}der paths to $\left(  2n+j,j\right)
$. Hence they satisfy the recursion%
\begin{equation}
D\left(  n,n+j\right)  =\omega D\left(  n-1,n-1+j\right)  +D\left(
n,n+j-1\right)  +D\left(  n-1,n+j\right) \label{(RecDellanoy)}%
\end{equation}

\[%
\begin{array}
[c]{c}%
\begin{array}
[c]{ccccccccc}%
j\uparrow &  &  &  &  &  &  & 1 & 0\\
5 &  &  &  &  &  & 1 & 0 & 7+6\omega\\
4 &  &  &  &  & 1 & 0 & 6+5\omega & 0\\
3 &  &  &  & 1 & 0 & 5+4\omega & 0 & 21+30\omega+10\omega^{2}\\
2 &  &  & 1 & 0 & 4+3\omega & 0 & 15+20\omega+6\omega^{2} & 0\\
1 &  & 1 & 0 & 3+2\omega & 0 & 10+12\omega+3\omega^{2} & 0 & 129\\
0 & 1 & 0 & 2+\omega & 0 & 6+6\omega+\omega^{2} & 0 & 63 & 0\\
n\rightarrow & 0 & 1 & 2 & 3 & 4 & 5 & 6 & 7
\end{array}
\\
\text{Uncompressed Grand Schr\"{o}der numbers (}\omega=1\text{)}%
\end{array}
\]

The generating function%
\[
\sum_{n=0}^{\infty}\sum_{l=0}^{n}\binom{k}{l}\binom{n+k-l}{n-l}\omega^{l}%
t^{n}=\frac{1}{1-t}\left(  \frac{1+t\omega}{1-t}\right)  ^{k}%
\]
shows that $D\left(  n,k\right)  $ is a Sheffer polynomial of degree $n$ in
$k$. The \emph{Delannoy polynomial} is of the form%
\begin{align*}
d_{k}\left(  t\right)   & =\sum_{j=0}^{k}D\left(  k-j,j\right)  t^{j}%
=\sum_{j=0}^{k}\sum_{l=0}^{k-j}\binom{j}{l}\binom{k-l}{j}\omega^{l}t^{j}\\
& =\sum_{l=0}^{k}\binom{k-l}{l}\omega^{l}t^{l}\sum_{j=0}^{k-2l}\binom{k-2l}%
{j}t^{j}=\sum_{l=0}^{k}\binom{k-l}{l}\omega^{l}t^{l}\left(  1+t\right)
^{k-2l}%
\end{align*}
and has generating function%
\[
\sum_{k=0}^{\infty}d_{k}\left(  t\right)  x^{k}=\frac{1}{1-x-t\left(  x+\omega
x^{2}\right)  }.
\]
From (\ref{(RecDellanoy)}) follows for $\omega=1$ that%
\[
d_{k-1}\left(  t\right)  =td_{k-1}\left(  t\right)  +td_{k-2}\left(  t\right)
+d_{k}\left(  t\right)
\]

Also for $\omega=1$ holds%

\[
p_{k}\left(  t\right)  =t^{-k}\sum_{l=0}^{k}\tbinom{k-l}{l}t^{l}\left(
1-t\right)  ^{k-2l}\left(  -1\right)  ^{l}=t^{-k}d_{k}\left(  -t\right)
\]
(see (\ref{(cp-p_n)})). Hence%
\begin{equation}
\mathcal{S}^{\left(  k\right)  }\left(  t;1\right)  =\sum_{n\geq0}%
S_{n;1}^{\left(  k\right)  }t^{n}=\frac{p_{k-1}\left(  t\right)  }%
{tp_{k}\left(  t\right)  }=\frac{d_{k-1}\left(  -t\right)  }{d_{k}\left(
-t\right)  }\label{(SchroederDelannoy)}%
\end{equation}
by (\ref{(GenFunS_n^(k))}). This shows an intimate connection between the
generating function of the Schr\"{o}der numbers in a band and the Delannoy
polynomials, when $\omega=1$. The Delannoy polynomials at negative argument,
$d_{k}\left(  -t\right)  $, satisfy for $\omega=1$ the same recursion as
$d_{k}\left(  t\right)  $,%
\[
d_{k-1}\left(  -t\right)  =td_{k-1}\left(  -t\right)  +td_{k-2}\left(
-t\right)  +d_{k}\left(  -t\right)  .
\]
This follows again from (\ref{(RecDellanoy)}).

Another connection exists with the inverse polynomial $s_{n}\left(  t\right)
$; from (\ref{(inverseSchroeder)}) $s_{n}\left(  t\right)  =\sum_{k\geq
0}s_{n,k}t^{n-k}=$%

\[
t\omega\sum_{m=1}^{n}\binom{n-m}{m-1}\left(  -1\right)  ^{m+1}\left(  1-\omega
t\right)  ^{n-2m}t^{m}-\sum_{m=0}^{n}\binom{n-m+1}{m}\left(  -1\right)
^{m+1}\left(  1-\omega t\right)  ^{n-2m}t^{m}%
\]
follows for $\omega=1$%

\begin{equation}
s_{n}\left(  t\right)  =\frac{t^{2}}{1-t}d_{n-1}\left(  -t\right)
+d_{n+1}\left(  -t\right)  /\left(  1-t\right)  .\label{(InverseDelannoy)}%
\end{equation}
and vice-versa,%

\begin{align*}
d_{n+1}\left(  -t\right)   & =\left(  1-t\right)  s_{n}\left(  t\right)
-t^{2}d_{n-1}\left(  -t\right) \\
& =\left(  1-t\right)  \sum_{i=0}^{n/2}t^{2i}\left(  -1\right)  ^{i}%
s_{n-2i}\left(  t\right)  +\left(  n\operatorname{mod}2\right)  \left(
-1\right)  ^{\left(  n+1\right)  /2}t^{n+1}%
\end{align*}
Hence the generating function of the bounded Schr\"{o}der numbers can for
$\omega=1$ be written as%
\[
\mathcal{S}^{\left(  k\right)  }\left(  t;1\right)  =\sum_{n\geq0}%
S_{n;1}^{\left(  k\right)  }t^{n}=\frac{\left(  1-t\right)  \sum
_{i=0}^{\left(  k-2\right)  /2}t^{2i}\left(  -1\right)  ^{i}s_{k-2-2i}\left(
t\right)  +\left(  k\operatorname{mod}2\right)  \left(  -1\right)  ^{\left(
k-1\right)  /2}t^{k-1}}{\left(  1-t\right)  \sum_{i=0}^{\left(  k-1\right)
/2}t^{2i}\left(  -1\right)  ^{i}s_{k-1-2i}\left(  t\right)  +\left(  \left(
k-1\right)  \operatorname{mod}2\right)  \left(  -1\right)  ^{k/2}t^{k}}%
\]

\section{Schr\"{o}der in a Band}

From (\ref{(SchroederDelannoy)}) follows for $\omega=1$ the generating
function of the (compressed) bounded (by $k$) Schr\"{o}der numbers,%

\begin{equation}
\mathcal{S}^{\left(  k\right)  }\left(  t;1\right)  =\frac{d_{k-1}\left(
-t\right)  }{d_{k}\left(  -t\right)  }\label{(genFuncBoundedSchroeder)}%
\end{equation}

Example: $\mathcal{S}^{\left(  4\right)  }\left(  t;1\right)  =\frac
{d_{3}\left(  -t\right)  }{d_{4}\left(  -t\right)  }=\allowbreak
\frac{1-5t+5t^{2}-t^{3}}{1-7t+13t^{2}-7t^{3}+t^{4}}$\newline$=\allowbreak
1+2t+6t^{2}+22t^{3}+89t^{4}+377t^{5}+1630t^{6}+\allowbreak7110t^{7}%
+31\,130t^{8}+136\,513t^{9}+599\,041t^{10}+\allowbreak2629\,418t^{11}%
+11\,542\,854t^{12}+50\,674\,318\allowbreak t^{13}+222\,470\,009t^{14}%
+976\,694\,489t^{15}+4287\,928\,678\allowbreak t^{16}+O\left(  t^{17}\right)
\allowbreak$%
\[%
\begin{tabular}
[c]{c}%
\begin{tabular}
[c]{l||lllllllll}\cline{2-10}%
$\uparrow m$ &  & \multicolumn{1}{|l}{} & \multicolumn{1}{|l}{} &
\multicolumn{1}{|l}{} & \multicolumn{1}{|l}{} & \multicolumn{1}{|l}{} &
\multicolumn{1}{|l}{} &  & \\\cline{2-10}%
k=4 &  & \multicolumn{1}{|l}{} & \multicolumn{1}{|l}{} & \multicolumn{1}{|l}{}
& \multicolumn{1}{|l}{} & \multicolumn{1}{|l}{} & \multicolumn{1}{|l}{} &  &
\\\cline{2-10}%
$3$ &  & \multicolumn{1}{|l}{} & \multicolumn{1}{|l}{} &
\multicolumn{1}{|l}{$1$} & \multicolumn{1}{|l}{$7$} & \multicolumn{1}{|l}{$36$%
} & \multicolumn{1}{|l}{$168$} & $756$ & $3353$\\\cline{2-10}%
$2$ &  & \multicolumn{1}{|l}{} & \multicolumn{1}{|l}{$1$} &
\multicolumn{1}{|l}{$6$} & \multicolumn{1}{|l}{$29$} &
\multicolumn{1}{|l}{$132$} & \multicolumn{1}{|l}{$588$} & $2597$ &
$11\,430$\\\cline{2-10}%
$1$ &  & \multicolumn{1}{|l}{$1$} & \multicolumn{1}{|l}{$4$} &
\multicolumn{1}{|l}{$16$} & \multicolumn{1}{|l}{$67$} &
\multicolumn{1}{|l}{$288$} & \multicolumn{1}{|l}{$1253$} & $5480$ &
$24\,020$\\\cline{2-10}%
$0$ & $1$ & \multicolumn{1}{|l}{$2$} & \multicolumn{1}{|l}{$6$} &
\multicolumn{1}{|l}{$22$} & \multicolumn{1}{|l}{$89$} &
\multicolumn{1}{|l}{$377$} & \multicolumn{1}{|l}{$1630$} & $7110$ &
$31130$\\\hline\hline
& $0$ & $1$ & $2$ & $3$ & $4$ & $5$ & $6$ & $7$ & $\rightarrow n$%
\end{tabular}
\\
\multicolumn{1}{l}{The compressed bounded ($k=4$) Schr\"{o}der numbers
($\omega=1$)}%
\end{tabular}
\]

From the recursion (\ref{(RecDellanoy)}) follows that
\begin{align*}
d_{k-1}\left(  -t\right)   & =td_{k-1}\left(  -t\right)  +td_{k-2}\left(
-t\right)  +d_{k}\left(  -t\right)  \text{ and therefor}\\
d_{k-1}\left(  -t\right)   & =\frac{t}{1-t}d_{k-2}\left(  -t\right)  +\frac
{1}{1-t}d_{k}\left(  -t\right)
\end{align*}

Hence
\[
d_{k-1}\left(  -t\right)  -td_{k-2}\left(  -t\right)  =\frac{t^{2}}%
{1-t}d_{k-2}\left(  -t\right)  +\frac{1}{1-t}d_{k}\left(  -t\right)
=s_{k-1}\left(  t\right)
\]
(see (\ref{(InverseDelannoy)})) and%
\[
\frac{d_{k-1}\left(  -t\right)  }{d_{k}\left(  -t\right)  }\left(
d_{k-1}\left(  -t\right)  -td_{k-2}\left(  -t\right)  \right)  -\frac
{d_{k-1}\left(  -t\right)  }{d_{k}\left(  -t\right)  }s_{k-1}\left(  t\right)
=0
\]
Therefore%
\[
\frac{d_{k-1}\left(  -t\right)  }{d_{k}\left(  -t\right)  }s_{k-1}\left(
t\right)  =\mathcal{S}^{\left(  k\right)  }\left(  t,k-1;1\right)  -\left(
d_{k-2}(-t)-td_{k-3}\left(  -t\right)  \right)
\]
(see (\ref{(genFuncBoundedSchroeder)}).

\begin{theorem}
\label{ThSchroeder}The power series part of $t^{-k}\mathcal{S}^{\left(
k\right)  }\left(  t;1\right)  s_{k-1}\left(  t\right)  $ equals
$t^{-k}S^{\left(  k\right)  }\left(  t,k-1;1\right)  $.
\end{theorem}

Example: (a) $t^{-4}\mathcal{S}^{\left(  4\right)  }\left(  t;1\right)
s_{3}\left(  t\right)  =\allowbreak\frac{\left(  t-1\right)  \left(
2t^{3}-8t^{2}+6t-1\right)  \left(  1-4t+t^{2}\right)  }{\left(  1-7t+13t^{2}%
-7t^{3}+t^{4}\right)  t^{4}}\allowbreak$\newline$=\allowbreak(t^{-4}%
-4t^{-3}+2t^{-2})+1+7t+36t^{2}+168t^{3}+\allowbreak756t^{4}+3353t^{5}%
+14\,783t^{6}+65\,016t^{7}+285\,648\allowbreak t^{8}+1254\,456t^{9}$%
\newline$+5508\,097t^{10}+24\,183\,271\allowbreak t^{11}+106\,173\,180t^{12}%
+O\left(  t^{13}\right)  \allowbreak$

(b) $t^{-4}\mathcal{S}^{\left(  4\right)  }\left(  t,4-1;1\right)
=\allowbreak t^{-4}\frac{1-5t+5t^{2}-t^{3}}{1-7t+13t^{2}-7t^{3}+t^{4}}\left(
1-6t+8t^{2}-2t^{3}\right)  -\frac{2t^{2}+1-4t}{t^{4}}=\allowbreak$%
\newline$\allowbreak1+7t+36t^{2}+168t^{3}+756t^{4}+3353t^{5}%
+14\,783\allowbreak t^{6}+65\,016t^{7}+285\,648t^{8}+1254\,456t^{9}$

$+5508\,097\allowbreak t^{10}+24\,183\,271t^{11}+106\,173\,180t^{12}+O\left(
t^{13}\right)  $

\end{document}